\documentclass{amsart}
\usepackage[latin1]{inputenc}
\usepackage{amssymb}

\newtheorem{theo}{Theorem}
\newtheorem{prop}[theo]{Proposition}
\newtheorem{lemma}[theo]{Lemma}

\def\NN{\mathbb N}
\def\ZZ{\mathbb Z}

\def\EE{\mathbb{E}}
\def\RR{\mathbb{R}}

\def\A{\mathcal{A}}
\def\U{\mathcal{U}}
\def\V{\mathcal{V}}
\def\W{\mathcal{W}}
\def\M{\mathcal{M}}
\def\N{\mathcal{N}}
\def\C{\mathcal{C}}
\def\P{\mathcal{P}}

\def\F{\mathcal{F}}

\def\B{\mathcal{B}}
\def\D{\mathcal{D}}
\def\uno{{ 1\!\!1}}
\def\funcion#1#2#3{{#1}\!:\!{#2}\rightarrow {#3}}
\def\'#1{{\if #1i{\accent"13\i}\else {\accent"13 #1}\fi}}
\def\conjunto#1#2{\{{#1}/{#2}\}}

\begin{document}
\title{Local conditional entropy in measure for covers with respect to a fixed partition}
\author{Pierre-Paul Romagnoli}
\address{Departamento de Matem\'aticas, Universidad Andres Bello, Rep\'ublica 252, Santiago, Chile}
\email{promagnoli@unab.cl}
\maketitle

\begin{abstract}
In this paper, we introduce two measure theoretical notions of
conditional entropy for finite measurable covers conditioned to a finite measurable partition and prove that they are equal.
Using this we state a local variational principle with respect to the notion of conditional entropy defined by Misiurewicz in \cite{M} for the case of open covers. This in particular extends the work done in \cite{R} and \cite{HYZ}.
\end{abstract}

\section{Introduction}

Topological dynamics and Ergodic Theory exhibit a remarkable parallelism. It is usual to find counterparts in both theories as is the case of transitivity, weak mixing and strong mixing and in particular with the notions that we study in this paper, this is topological entropy and measure theoretical entropy. Although, even if the notions and results are similar the methods to prove them are quite different. In many cases the connection between these notions is through a variational principle by taking suprema over all invariant measures to obtain the topological notion. For the case of $\ZZ$ actions  Measure-theoretical Entropy for an invariant measure was introduced in 1958 in \cite{K} and
Topological Entropy in 1965 in \cite{AKM}. In 1969 and 1970 Goodman \cite{G} and Goodwyn \cite{Gw} in two separate papers proved the first variational principle. M. Misiurewicz introduced the notion of topological conditional entropy in \cite{M} and the variational principle was established by F. Ledrappier in \cite{L}.\\

This paper is inserted in the so called local theory of entropy for Topological Dynamical Systems that started in the early 90's with the work of Fran\c cois Blanchard (see \cite{B} and \cite{BL}) that is interesting by itself but has also proven to be fundamental to many other related areas. For instance the existence of topological Pinsker factors, zero entropy factors, disjointness theorems, characterizations of positive entropy, entropy pairs and tuples and so on.\\

In \cite{BGH} while extending the notion of topological entropy pairs to a measure theoretical setting the authors proved a \textit{variational principle
for open covers} for a topological dynamical system $(X,T)$. More
precisely, for every open cover there exists a $T$-invariant
measure $\mu$ such that the topological entropy of the cover is
bounded above by the $\mu$-entropy of every partition finer than
the cover.

In \cite{R} the author gave a new approach to the variational principle given in \cite{BGH} extending the measure theoretical notions from partitions to covers and thus proposing two measure theoretical notions of
entropy for open covers denoted $h_\mu^+$ and $h_\mu^-$. In doing so he was then able to state and prove a local variational
principle for $h_\mu^-$. Namely, for every open cover
there exists a $T$-invariant measure $\mu$ such that the
topological entropy is equal to the $h_\mu^-$-entropy of the
cover. The notion $h_\mu^-$ verifies most of the relevant
properties of the usual notion of entropy for measurable
partitions. It is also proven that the infimum that appears in the
definition of $h_\mu^-$ is attained. This fact implies some good
properties, for instance it implies that $h_\mu^-$ is preserved by
measure-theoretic extension and as a function of $\mu$, $h_\mu^-$
is upper semicontinuos (see \cite{HY}). The main open question
left in \cite{R} was to know if the two notions were equal.

In \cite{HMRY} the authors proved that these notions were positive
simultaneously and, using the Jewett-Krieger Theorem they proved
that the equality between $h_\mu^+$ and $h_\mu^-$ is equivalent to
the existence of a local variational principle for $h_\mu^+$
(similar to the one proven for $h_\mu^-$). In \cite{GW} such a
variational principle was established for $h_\mu^+$ and thus the
equality holds. An alternative definition for measure entropy for covers was given in \cite{S} giving a new proof of the variational 
principle and at the same time proving the equality of $h_\mu^+$ and $h_\mu^-$ once again.

In \cite{HY} the authors applied the same ideas to extend the notion of topological pressure to a measure theoretical setting proving again a local variational principle of the
same topological pressure and in \cite{HYZ} a similar one in a conditional version
with respect to a factor. 

For the existence of a variational principle the requirement of the cover to be open is unavoidable. However the equality between $h_\mu^+$ and $h_\mu^-$ can be stated in general for measurable covers extending the topological notion to this kind of covers since it is turns out to be mainly combinatorial. This was well known even in the first known proofs of the global variational principle. The global bound of the variational principle fails even for nice non open covers as is the case of closed covers where for example the  existence of one non recurrent point implies that the suprema of the entropy over all closed covers is infinite. This was reviewed and explained by Goodwyn in \cite{Gw2} (see example 3) although it was already mentioned by him in \cite{Gw}.

Since then, the research has mainly focused in extending the notions of entropy and pressure to more general actions. In \cite{HYZ2} it was done for the case countable discrete amenable groups by using the tool of F{\o}lner sequences making the techniques and the proofs very similar. The case of continuous bundle
random dynamical systems of an infinite countable discrete amenable group action has been proven in \cite{DZ} and for sofic group actions in \cite{Z}.

However the main idea that gives the natural bridge between the measurable and the topological notions in every one of the versions of variational principles discussed so far is the one given in \cite{R}. This is, extended the measurable formulas that apply to partitions to covers by considering the infima over all partitions finer than the cover. Then for measurable covers the natural ``+" and ``-" definitions prove that they are equal for all measurable covers and then use this to prove a variational principle for open covers for the suitable topological equivalent.\\

Unfortunately the real local variational problem is still open, this is, with respect to a fixed open cover conditioned to another fixed open cover. The missing ingredient is a clear way to extend the notion with respect to the conditioning partition. Extending the definition with respect to the conditioning variable in the same way by taking infima is not the way to go for a local definition and even for a fixed partition using the conditional measurable definition of partitions is not enough. Considering the alternative definition of conditional entropy as an average of the entropy of the first partition with respect to the induced measures over the atoms of the conditioning partition and using infima in a clever way solves the issue giving once again two different definitions that turn out to be the same. In this paper we address this for the first time proving the local conditional version with respect to a finite measurable partition.\\

More precisely, in this paper we propose two notions of conditional measure
theoretical entropy for measurable finite covers conditioned to a fixed measurable partition that extends the notions
$h_\mu^-$ and $h_\mu^+$ given in \cite{R} and \cite{HYZ} and proving in general that they coincide for every finite measurable cover. In the case of finite open covers we prove that they satisfy a local variational principle with the notion of
conditional entropy defined by Misiurewicz in \cite{M} once again extending the results in \cite{R} and \cite{HYZ}.\\

\section{Basic definitions and results}

We will define a series of notions as limits that are also infima using the classic subadditive lemma that we state without proof:

\begin{lemma}
\label{subadd}
For any subadditive sequence $\{a_N\}_{N\in \NN}\subseteq \RR^+$, this is $a_{N+M}\le a_N+a_M$ for any $N,M\in \NN$ we have:

  \begin{equation}
  \lim\limits_{N\to\infty}{1\over N}a_N=\inf\limits_{N\in\NN}{1\over N}a_N.
  \end{equation}
   \end{lemma}
{\bf Proof:}
See \cite{D}.
\qed

Let us introduce the basic notation and definitions used in this
article. For more details on measure-theoretical and topological
entropy including the most recent results of the beginning of the year 2000
we refer to the deeply insightful and inspirational textbook by Tomasz Downarowicz \cite{D}.\\

 Let $(X,T)$ be a
topological dynamical system (TDS). This is, $X$ a compact metric
space and $\funcion{T}{X}{X}$ a homeomorphism. A TDS $(X,T)$ is
0-dimensional when the space $X$ has a countable open-closed (clopen)
topological basis. The set $\mathcal{M}_T(X)$ of $T$-invariant Borel
probability measures is a convex, compact (in the weak topology) and nonempty set.
We denote $\mathcal{M}^e_T(X)$ the set of ergodic measures. A measure theoretic dynamical system (MDS), $(X,\B,\mu,T)$, is a
probability space $(X,\B,\mu)$ and a bi-measurable bijection
$T:X\rightarrow X$ that preserves the measure $\mu$. So a TDS $(X,T)$ gives a family of MDS indexed by the set $\mathcal{M}_T(X)$ where $\B$ is the $\sigma$-algebra of the borel sets of $X$.\\

In this article a cover of $X$ is a finite
cover of Borel subsets of $X$. A cover is said to be an \emph{open cover}
if it consists only of open sets. A \emph{partition of $X$} is a cover by
pairwise disjoint sets. Let $\mathcal{P}_X$ denote the set of
partitions of $X$, $\mathcal{C}_X$ the set of covers of $X$ and
$\mathcal{C}^o_X$ the set of open covers of $X$.\\

 Given two covers
$\mathcal{U},\mathcal{V}\in \mathcal{C}_X$, $\mathcal{U}$ is said
to be finer than $\mathcal{V}$ ($\mathcal{U}\succeq \mathcal{V}$)
if for every $U\in \mathcal{U}$, there is $V\in \mathcal{V}$ such
that $U\subseteq V$. Let
$\mathcal{U}\vee\mathcal{V}=\conjunto{U\cap V}{U\in
\mathcal{U},V\in \mathcal{V}}$. It is clear that
$\mathcal{U}\vee\mathcal{V}\succeq \mathcal{U}$ and
$\mathcal{U}\vee\mathcal{V}\succeq \mathcal{V}$. However,
$\mathcal{V}\succeq \mathcal{U}$ does not imply that
$\mathcal{U}\vee\mathcal{V}=\mathcal{ V}$. Given integers $M\le N$
and $\mathcal{U}\in\mathcal{C}_X$, one sets,
$\mathcal{U}_M^{N}=\bigvee\limits_{n=M}^{N}T^{-n}\mathcal{U}$. For $\U=\{U_1,\dots,U_M\}$ define $\U^*=\conjunto{\alpha=\{A_1,\dots,A_M\}\in \P_X}{\forall m\in[M],A_m\subseteq U_m}$.\\

The set $\U^*$ is useful to simplify several proofs by using lemma 2 in \cite{HMRY} that we state without proof.\\

\begin{lemma} For any measure space $(X,\B,\mu)$ and $\funcion{G}{\P_X}{\RR}$ such that $\forall \alpha,\beta \in \P_X$, $\alpha\succeq\beta$ implies that $G(\alpha)\ge G(\beta)$ one has that for any $\U\in\C_X$:

\begin{equation}
\label{G}
\inf\limits_{\alpha\succeq \U\atop \alpha\in\P_X}G(\alpha)=\inf\limits_{\alpha\in\U^*}G(\alpha).
\end{equation}

\end{lemma}

In general one can define the $\mu$-conditional entropy of
$\alpha\in\mathcal{P}_X$ conditioned
by a sub $\sigma$-algebra $\A\subseteq \B$ as:

\begin{equation}
\label{defhmucondsigma}
H_\mu(\alpha|\A):=\sum\limits_{A\in \alpha}\int\limits_{X}\phi[\EE_\mu(\uno_A|\A)(x)]\,d\mu(x).
\end{equation}

Where $\funcion{\phi}{[0,1]}{\RR^+}$ is defined as $\phi(x):=\-x\log x$ for $x>0$ and $\phi(0)=0$. \\

Given
$\alpha,\beta\in\mathcal{P}_X$ the entropy of $\alpha$ conditioned
by $\beta$ is given by $H_\mu(\alpha|\beta)=H_\mu(\alpha|\overline{\beta})$ where $\overline{\beta}$ is the $\sigma$-algebra generated by $\beta$.
In this case there are simpler formulas to compute it, for instance when $\beta=\{X\}$ this gives the
entropy of $\alpha$ and $H_\mu(\alpha):=-\sum\limits_{A\in \alpha}\phi(\mu(A))$. In general for any $\beta\in\P_X$:

\begin{equation}
\label{defhmucond}
H_\mu(\alpha|\beta):=H_\mu(\alpha\vee\beta)-H_\mu(\beta)=
\sum_{B\in \beta}\mu(B)H_{\mu_B}(\alpha).
\end{equation}

where $\mu_B$ denotes the conditional measure induced by $\mu$ on
$B$ (zero if $\mu(B)=0$). The function $H_{\cdot}(\alpha|\A)$ is concave on $\M_T(X)$
(since $-x\log x$ is concave). \\

{\bf REMARK:} The formula $H_\mu(\alpha)$ can be applied to any disjoint family of measurable sets even if it does not cover $X$.\\

Given $\alpha,\beta\in \P_X$ or any sub-$\sigma$-algebra $A\subseteq \B$ by lemma \ref{subadd} the $\mu$-entropy of $\alpha$ conditioned to the partition $\beta\in
\P_X$ (or some $T$-invariant $\sigma$-algebra $\A$) with respect to $T$ is well defined as:

\begin{equation}
\label{defhmucondpart}
\begin{split}
h_\mu(\alpha|\beta,T)&:=\lim_{N\rightarrow\infty}\frac{1}{N}H_\mu(\alpha_0^{N-1}|\beta_0^{N-1})=
\inf_{N\in \NN}\frac{1}{N} H_\mu(\alpha_0^{N-1}|\beta_0^{N-1}),\\
h_\mu(\alpha|\A,T)&:=\lim_{N\rightarrow\infty}\frac{1}{N}H_\mu(\alpha_0^{N-1}|\A)=
\inf_{N\in \NN}\frac{1}{N} H_\mu(\alpha_0^{N-1}|\A).
\end{split}
\end{equation}

When $\beta=\{X\}$ this yields the standard $\mu$-entropy of $\alpha$ with respect to $T$ that we denote as
$h_\mu(\alpha,T):=h_\mu(\alpha|\{X\},T)$.\\

For any $\U=\{U_1,\dots U_M\},\V=\{V_1,\dots V_M\}\in \C_X$ define $\mu(\U\Delta\V)=\sum\limits_{m=1}^M \mu(U_m\Delta V_m)$. This definition gives a notion of distance between covers and partitions that is compatible with the conditional measure entropy as the following lemma whose proof is taken from Peter Walter's classic textbook \cite{W}.\\

\begin{lemma}
\label{conthcond}
Fix $M\in\NN$ and $\epsilon>0$ then there exists $\delta>0$ such that $\forall \alpha,\alpha',\beta\in \P_X$ with $|\alpha|=|\alpha'|=M$ if $\mu(\alpha\Delta\alpha')<\delta$ one has that $|h_\mu(\alpha|\beta,T)-h_\mu(\alpha'|\beta,T)|\le \epsilon$.
\end{lemma}
{\bf Proof:} See \cite{W}.
\qed

We use the same definition of chapter 6.3 in \cite{D} but for $\C_X$ and $\beta\in \P_X$:

\begin{equation}
\label{condnu1}
\mathcal{N}(\mathcal{U}|\beta):=
\max\limits_{B\in\beta}\mathcal{N}(\mathcal{U}\cap B),\>
\mathcal{N}(\mathcal{U}\cap B):=\min\conjunto{|\V|}{\V\subseteq \U, B\subseteq \bigcup\limits_{V\in \V}V}.
\end{equation}

From the same ideas used in Fact 6.3.2 in \cite{D} and other simple calculations:\\

\begin{prop}\label{propNcond}
Let $(X,T)$ be a TDS and $\U,\V\in \C_X$ and $\alpha,\beta,\gamma\in \P_X$ then\\

\begin{enumerate}
\item $\N(\U|\beta)=0 $ iff $ \beta\succeq \U$.
\label{propNcond1}

\item $\N(T^{-1}\U| T^{-1}\beta)=\N(\U|\beta)$.
\label{propNcond4}

\item $\N(\U|\beta)\ge \N(\V|\beta)\quad\rm{if}\quad\U\succeq\V
$.
\label{propNcond2}

\item $\N(\U|\beta)\leq
\N(\U|\gamma)\quad\rm{for}\quad\beta\succeq\gamma$.
\label{propNcond3}


\item $\N(\U\vee\V|\beta)\leq \N(\U|\beta)\cdot \N(\V|\beta)$.
\label{propNcond6}



\end{enumerate}

\end{prop}

{\bf Proof:} Part (\ref{propNcond4}) is trivial and for all the others use Fact 6.3.2. with $F=X$ and consider $\V$ and $\W$ as partitions when required.
\qed

Proposition
\ref{propNcond} parts (\ref{propNcond4}),(\ref{propNcond3}) and
(\ref{propNcond6}) imply that the sequence
$\{\log\N(\U_0^{N-1}|\beta_0^{N-1})\}_{N\in\NN}$ is subadditive, and
using lemma \ref{subadd} as in Definition 6.3.14 in \cite{D} we obtain the combinatorial (topological if $\U$ is open) entropy of the cover $\U\in \C_X$ conditioned by partition $\beta\in\P_X$ with respect to $T$ as:

\begin{equation}
\label{defhtopcond}
h(\U|\beta,T):=\lim_{N\rightarrow\infty}\frac{1}{N}\log \mathcal{N}(\U_0^{N-1}|\beta_0^{N-1})=
\inf_{N\in \NN}\frac{1}{N} \log
\mathcal{N}(\U_0^{N-1}|\beta_0^{N-1}).
\end{equation}

 When $\beta=\{X\}$ we recover the classical topological entropy and denote $\mathcal{N}(\U|\{X\})=\mathcal{N}(\U)$ and
$h(\U|\{X\},T)=h(\U,T)$.\\

In \cite{R} we defined for $\U\in \C^o_X$ and $\mu\in \M_T(X)$, $H_\mu(\U):=\inf\limits_{\alpha\succeq \U\atop \alpha\in \P_X}H_\mu(\alpha)$. To extend this idea we define for any $\U\in \C_X$ and $\beta\in\P_X$ given $B\in \beta$,  $H_{\mu_B}(\U)$ as for $H_\mu(\U)$ (even if $\mu_B$ is not $T$-invariant) and using the decomposition shown in (\ref{defhmucond}):

\begin{equation}
\label{defhmucondrec}
H_\mu(\U|\beta):=\sum\limits_{B\in \beta}\mu(B)H_{\mu_B}(\U),\> H^*_\mu(\U|\beta):=\inf\limits_{\alpha\succeq \U}H_{\mu}(\alpha|\beta).
\end{equation}

It is straightforward that $H_\mu(\U|\beta)\le H^*_\mu(\U|\beta)$ and when $\beta=\{X\}$ and $\U\in \C_X^o$ we recover in both cases the definition given in \cite{R} since $H_\mu^*(\U|\{X\})=H_\mu(\U|\{X\})=H_\mu(\U)$. Nevertheless we will prove that they are actually always equal. This result and some basic properties for the conditional entropy will extend to $H_\mu(\U|\beta)$ with respect to any fixed $\beta\in \C_X$ with the use of two key facts. From \cite{R} as mentioned in Fact 8.3.5 in \cite{D},

\begin{equation}
  \label{extpart}
H_{\mu}(\U)=\min\limits_{\alpha\in Ext(\U)}H_{\mu_B}(\alpha).
\end{equation}

With
 $Ext(\U)$ the set of partitions of the form $\{U_1,U_2\backslash U_1,\dots,U_d\backslash (U_1\cup\dots\cup U_{d-1})\}$ where
 $\{U_1,\dots,U_d\}$ is an ordering of $\U$. This is also true for any $B\in\beta$ for $H_{\mu_B}(\U)$ since the proof does not need $\mu$ to be $T$-invariant.\\

 When $A$ and $B$ are disjoint then:
\begin{equation}
\label{aubhmu}
\mu(A\cup B)H_{\mu_{A\cup B}}(\alpha)\ge \mu(A)H_{\mu_{A}}(\alpha)+\mu(B)H_{\mu_{B}}(\alpha)\ge \mu(A)H_{\mu_{A}}(\alpha).
\end{equation}

Thus if $A\subseteq B$ then $\mu(A)H_{\mu_{A}}(\U)\le \mu(B)H_{\mu_{B}}(\U)$.\\

\begin{lemma}
For any $\U\in \C_X$, $\beta\in \P_X$, $H_\mu(\U|\beta)=H^*_\mu(\U|\beta)$.
\end{lemma}

{\bf Proof:}
Fix $\U\in \C_X$ and $\beta\in\P_X$, we only need to prove that $H^*_\mu(\U|\beta)\le H_\mu(\U|\beta)$. From fact \ref{extpart} for any $B\in \beta$ there exists $\alpha_B\in \P_X$ such that $H_{\mu_B}(\U)=H_{\mu_B}(\alpha_B)$ so $H_\mu^*(\U|\beta)=\sum\limits_{B\in\beta}\mu(B)H_{\mu_B}(\alpha_B)$.\\

Define $\alpha=\bigvee\limits_{B\in \beta}[\alpha_B\cap B]\cup \{B^c\}$. By definition for any any $\tilde\alpha\succeq \U$ and $B\in \beta$, $H_{\mu_B}(\tilde \alpha)\ge H_{\mu_B}(\alpha_B)=H_{\mu_B}([\alpha_B\cap B]\cup \{B^c\})=H_{\mu_B}(\alpha)$ so multiplying by $\mu(B)$ and adding up over $B\in\beta$ we conclude that $H_\mu(\U|\beta)=H_\mu(\alpha|\beta)=H^*_\mu(\U|\beta)$.\\
\\






\qed


{\bf REMARK:} From now we will use as definition of $H_\mu(\U|\beta)$ the most suitable of these two formulas as needed in the proofs.\\

Now we prove that this notion has the standard properties of a static entropy.\\

\begin{prop}\label{propHmur}
Let $(X,\B,\mu,T)$ be a MDS and $\U,\V,\W\in \C_X$ and $\alpha,\beta,\gamma\in \P_X$ then\\

\begin{enumerate}
\item $0\le H_\mu(\U|\beta) \le \log \mathcal{N}(\U|\beta)$ and
$H_\mu(\U|\beta)=0 $ if $ \beta\succeq \U$.
\label{propHmur1}

\item $H_\mu(T^{-1}\U| T^{-1}\beta)= H_\mu(\U|\beta)$.
\label{propHmur4}


\item $H_\mu(\U|\beta)\leq
H_\mu(\U|\gamma)\quad\rm{for}\quad\beta\succeq\gamma$.
\label{propHmur3}


\item $H_\mu(\U\vee\V|\beta)\leq H_\mu(\U|\beta)+H_\mu(\V|\beta)$.
\label{propHmur6}



\end{enumerate}

\end{prop}

{\bf Proof:} Fix $\U,\V,\W\in \C_X$ and $\alpha,\beta,\gamma\in \P_X$.\\

$(\ref{desmutopcond})$ For any  $B\in \beta$,
$H_{\mu_B}(\alpha)\le \log|\alpha\cap B|$ and $\inf\limits_{\alpha\succeq \mathcal{U}}|\alpha\cap B|
=\mathcal{N}(\mathcal{U}\cap B)$. So:

\begin{equation}
H_\mu(\mathcal{U}|\beta)=\sum_{B\in
\beta}\mu(B)H_{\mu_B}(\U)\le\sum_{B\in
\beta}\mu(B)\log \mathcal{N}(\mathcal{U}\cap B)\le\log\mathcal{N}(\mathcal{U}|\beta).
\end{equation}

If $\beta\succeq\U$, $H_{\mu_B}(\U)\le H_{\mu_B}(\beta)=0$ for any $B\in \beta$.\\

$(\ref{propHmur4})$ For any  $B\in \beta$ from equation (\ref{extpart}):

\begin{align*}
H_{\mu_{T^{-1}B}}(T^{-1}\U)&=\min\limits_{\alpha\in Ext(T^{-1}\U)}H_{\mu_{T^{-1}B}}(\alpha)=\min\limits_{\alpha\in T^{-1}Ext(\U)}H_{\mu_{T^{-1}B}}(\alpha),\\
&=\min\limits_{\alpha\in Ext(\U)}H_{\mu_{T^{-1}B}}(T^{-1}\alpha)=\min\limits_{\alpha\in Ext(\U)}H_{\mu_{B}}(\alpha)=H_{\mu_B}(\U).
\end{align*}

Adding up over all $B\in\beta$ and since $\mu(T^{-1}B)=\mu(B)$ we conclude the result. \\

$(\ref{propHmur3})$ Let $\beta=\{B_1,\dots,B_N\}$ and $\gamma=\{C_1,\dots,C_M\}$ if $\beta\succeq \gamma$ then there exists a partition $\{K_n\}_{n=1}^N$ of $\{1,\dots,N\}$ such that $B_n=\bigcup\limits_{k\in K_n}C_k$ for any $1\le n\le N$.\\

From equation (\ref{aubhmu}) for any $1\le n\le N$:

\begin{equation}
\mu(B_n)H_{\mu_{B_n}}(\U)\ge \sum\limits_{k\in K_n}\mu(C_n)H_{\mu_{C_n}}(\U).
\end{equation}

Adding up over $n$ concludes the result.\\

$(\ref{propHmur6})$ For any $\alpha\succeq \U$ and $\gamma\succeq \V$ from basic properties of $H_\mu$ for partitions:

\begin{equation}
H_\mu(\U\vee\V|\beta)\le
H_\mu(\alpha\vee\gamma|\beta)\le
H_\mu(\alpha|\beta)+H_\mu(\gamma|\beta).
\end{equation}

Taking infima over all $\alpha\succeq \U$ and $\gamma\succeq \V$ concludes the result.\\

\qed

This allows us to define two measure theoretical notions of dynamical entropy for a measurable cover conditioned to a fixed partition.\\

Proposition
\ref{propHmur} parts (\ref{propHmur4}),(\ref{propHmur3}) and
(\ref{propHmur6}) imply that the sequence
$\{H_\mu(\U_0^{N-1}|\beta_0^{N-1})\}_{N\in\NN}$ is subadditive, and
so by lemma \ref{subadd} we define two notions of $\mu$-entropy of the cover $\U\in \C_X^o$ conditioned by the
partition $\beta\in \P_X$ with respect to $T$ as:

\begin{equation}
\label{defhmucondh}
\begin{split}
h^-_\mu(\U|\beta,T)&:=\lim_{N\rightarrow\infty}\frac{1}{N}H_\mu(\U_0^{N-1}|\beta_0^{N-1})=
\inf\limits_{N\in \NN}\frac{1}{N} H_\mu(\U_0^{N-1}|\beta_0^{N-1}), \\
h^+_\mu(\U|\beta,T)&:=\inf\limits_{\alpha\succeq \U\atop \alpha\in\P_X}h_\mu(\alpha|\beta,T).
\end{split}
\end{equation}

Notice that when $\beta=\{X\}$ we recover the notions in \cite{R}.\\


The following lemma is extremely useful to extend properties between $h_\mu^+$ and $h_\mu^-$ and finally prove they coincide and it will be used many times in the sequel.\\

\begin{lemma}
\label{lem:pvc-df>}
Let $(X,\mathcal{B},\mu,T)$ a MDS and $\mathcal{U}\in \mathcal{C}_X$ and $\beta\in \P_X$.
Then:\\

\begin{enumerate}
\item $h^-_\mu(\U|\beta,T)\le h^+_\mu(\U|\beta,T)$ and $h^-_\mu(\U|\beta,T)\le h(\U|\beta,T)$.
\label{desmutopcond}
\item $\forall M\in\NN$, $h^-_\mu(\U|\beta,T)= {1\over M}h^-_\mu(\U_0^{M-1}|\beta_0^{M-1},T^M)$.
\label{igmucondtm}
\item $h^-_\mu(\U|\beta,T)= \lim\limits_{M\to\infty}{1\over M}h^+_\mu(\U_0^{M-1}|\beta_0^{M-1},T^M)$.
\label{limtmhmcond}
\end{enumerate}

\end{lemma}

{\bf Proof:}
Fix $\U\in \C_X$ and $\beta\in\P_X$.\\

(\ref{desmutopcond}) By Proposition \ref{propHmur} part (\ref{propHmur1}), $H_\mu(\U|\beta) \le \log \mathcal{N}(\U|\beta)$ so:

\begin{align*}
h_\mu^-(\mathcal{U}|\beta,T)&=\lim_{N\rightarrow \infty}\frac{1}{N}H_\mu(\mathcal{U}_0^{N-1}|\beta_0^{N-1}),\\
&\le\lim_{N\rightarrow
\infty}\frac{1}{N}\log\mathcal{N}(\mathcal{U}_0^{N-1}|\beta_0^{N-1})=h(\mathcal{U}|\beta,T).
\end{align*}

For any $N\in\NN$, $H_\mu(\U_0^{N-1}|\beta_0^{N-1})\le H^*_\mu(\U_0^{N-1}|\beta_0^{N-1})\le \inf\limits_{\alpha\succeq \U}H_\mu(\alpha_0^{N-1}|\beta_0^{N-1})$. Dividing by $N$ and taking limit proves that $h^-_\mu(\U|\beta,T)\le h^+_\mu(\U|\beta,T)$. \\

(\ref{igmucondtm}) For every $M\in\NN$ just by definition:

\begin{align*}
h^-_\mu(\U|\beta,T)&=\lim\limits_{N\to\infty}{1\over NM}H_\mu(\U_0^{NM-1}|\beta_0^{NM-1}),\\
&={1\over M}\lim\limits_{N\to\infty}{1\over N}H_\mu\left(\bigvee\limits_{n=0}^{N-1}T^{-nM}\U_0^{M-1}|
\bigvee\limits_{n=0}^{N-1}T^{-nM}\beta_0^{M-1}\right),\\
&={1\over M}h^-_\mu(\U_0^{M-1}|\beta_0^{M-1},T^M).
\end{align*}

 (\ref{limtmhmcond}) From parts (\ref{desmutopcond}) and (\ref{igmucondtm}) for every $M\in\NN$:

\begin{equation}
h^-_\mu(\U|\beta,T)= {1\over M}h^-_\mu(\U_0^{M-1}|\beta_0^{M-1},T^M)\le {1\over M}h^+_\mu(\U_0^{M-1}|\beta_0^{M-1},T^M).
\end{equation}

 Also for every $N\in\NN$ using Proposition \ref{propHmur} parts (\ref{propHmur3}) and (\ref{propHmur6}):

 \begin{align*}
{1\over M}h^+_\mu(\U_0^{M-1}|\beta_0^{M-1},T^M)&=\lim\limits_{N\to\infty}{1\over NM}\inf\limits_{\alpha\succeq \U_0^{M-1}}H_\mu\left(\bigvee\limits_{n=0}^{N-1}T^{-nM}\alpha|\beta_0^{NM-1}\right),\\
&\le\lim\limits_{N\to\infty}{1\over NM}\inf\limits_{\alpha\succeq \U_0^{M-1}}\sum\limits_{n=0}^{N-1} H_\mu\left(T^{-nM}\alpha|T^{-nM}\beta_0^{M-1}\right),\\
&\le{1\over M}\inf\limits_{\alpha\succeq \U_0^{M-1}} H_\mu\left(\alpha|\beta_0^{M-1}\right)={1\over M}H_\mu\left(\U_0^{M-1}|\beta_0^{M-1}\right).
 \end{align*}

Taking limit as $M$ goes to infinity concludes the proof.

\qed

\begin{prop}
\label{hmufactors}
Let $(X, T)$ and $(Y,S)$ be TDS, $\mu\in\M_T(X)$ and
$\nu\in\M_S(Y)$.Let $\funcion{\varphi}{(X, T,\mu)}{(Y,S,\nu)}$ be a measure-theoretical factor map,
$\U\in \C_Y(S)$ and $\beta\in \P_Y$. Then $h^-_\nu(\U|\beta,S)=h^-_\mu(\varphi^{-1}\U|\varphi^{-1}\beta,T)$ and $h(\U|\beta,S)=h(\varphi^{-1}\U|\varphi^{-1}\beta,T)$.
\end{prop}

{\bf Proof:}
Fix $\U\in \C_Y(S)$ and $\beta\in \P_Y$. Clearly for any $N\in\NN$, $\N(\varphi^{-1}\U_0^{N-1}|\varphi^{-1}\beta_0^{N-1})=\N(\U_0^{N-1}|\beta_0^{N-1})$ so taking $\log$ diving by $N$ and taking limit as $N$ tends to infinity proves that $h(\U|\beta,S)=h(\varphi^{-1}\U|\varphi^{-1}\beta,T)$.\\


 In Proposition 6 in \cite{R} it is proven that $H_\mu(\varphi^{-1}\U)=H_\nu(\U)$ and that part of the proof does not require for $\mu$ to be $T$-invariant only that $\mu\circ \varphi^{-1}=\nu$. Since for any measurable set $B$ it is straightforward to see that $\mu_{\varphi^{-1}B}\circ \varphi^{-1}=\nu_{B}$ thus
the same proof applies to prove that $H_{\mu_{\varphi^{-1}B}}(\varphi^{-1}\U)=H_{\nu_B}(\U)$.
Taking the sum over all $B\in \beta$ we conclude that:

\begin{equation}
\label{hnuphi=hphimu}
H_\mu(\varphi^{-1}\U|\varphi^{-1}\beta)=H_\nu(\U|\beta).
\end{equation}

Using equation (\ref{hnuphi=hphimu}) since $S\circ \varphi=\varphi\circ T$ for any $N\in\NN$ we have that:

\begin{equation}
H_\mu\left(\bigvee\limits_{n=0}^{N-1}T^{-n}(\varphi^{-1}\U)|\bigvee\limits_{n=0}^{N-1}T^{-n}\varphi^{-1}\beta\right)=
H_\nu\left(\bigvee\limits_{n=0}^{N-1}S^{-n}\U|\bigvee\limits_{n=0}^{N-1}S^{-n}\beta\right).
\end{equation}

Dividing and taking the limit when $N$ goes to infinity concludes the proof.\\
\qed

\section{Local Variational Principles}

In this section we prove the local variational principle for both notions. The techniques and proofs are a mix of the work done in \cite{R},\cite{GW},\cite{HMRY} and \cite{HYZ}.\\

\begin{lemma}\label{lemmappio+}
Let $(X,T)$ a TDS, $\mathcal{U}\in \mathcal{C}_X$
and $\beta\in\mathcal{P}_X$. For every family
$\{\alpha_l\}_{l=1}^K$ of finite partitions finer than
$\mathcal{U}$, for every $N\in \NN$ choose $B\in \beta_0^{N-1}$
such that the maximum of
$\mathcal{N}(\mathcal{U}_0^{N-1}|\beta_0^{N-1})$ is attained.

There exists
a finite subset $P\subseteq B$ such that every element of $(\alpha_l)_0^{N-1}$
contains at most one point of $P$ for every $1\le l\le K$ and
$|P|\ge \frac{\mathcal{N}(\mathcal{U}_0^{N-1}|\beta_0^{N-1})}{K}$.
\end{lemma}

{\bf Proof:} Choose $\U\in \C_X$, $\beta\in \P_X$, $N\in\NN$ and $B\in \beta_0^{N-1}$ as stated.\\

For $x\in X$ and $1\le
l\le K$, let $A_l(x)$ be the element in $(\alpha_l)_0^{N-1}$ that
contains $x$.
Take $x_1\in B$. If $B\subseteq \bigcup_{l=1}^K A_l(x_1)$
then let $P=\{x_1\}$. Every element in $(\alpha_l)_0^{N-1}$
contains at most one point of $P$ for every $1\le l\le K$ and
$|P|K \ge \mathcal{N}(\mathcal{U}_0^{N-1}|\beta_0^{N-1})$.

If $B\nsubseteq
\bigcup_{l=1}^K A_l(x_1)$ let $B_1=B\setminus \bigcup_{l=1}^K
A_l(x_1)$ and take $x_2\in B_1$. If $B_1\subseteq \bigcup_{l=1}^K A_l(x_2)$
then $P=\{x_1,x_2\}$. Every element in $(\alpha_l)_0^{N-1}$
contains at most one point of $P$ for every $1\le l\le K$ and $|P_N|K \ge \mathcal{N}(\mathcal{U}_0^{N-1}|\beta_0^{N-1})$.\\

Otherwise,
$B_2=B_1\setminus \bigcup_{l=1}^K A_l(x_2).$
Since this is a finite procedure we obtain a set
$\{x_1,\ldots x_m\}$ such that $B\subseteq
\bigcup_{j=1}^{m}\bigcup_{l=1}^K
A_l(x_j)$.

Let $P=\{x_1,\ldots x_m\}$. By construction $(\alpha_l)_0^{N-1}$
contains at most one point of $P$ for every $1\le l\le K$  and
$|P|K=|\bigcup_{j=1}^{m}\bigcup_{l=1}^K A_l(x_j)|\ge
\mathcal{N}(\mathcal{U}_0^{N-1}|\beta_0^{N-1}).$
\qed



\begin{theo}
\label{varprinc}
For every TDS $(X,T)$, $\U\in \C^o_X$ and $\beta\in \P_X$ there exists $\mu\in\M_T(X)$ such that $h_\mu^+(\U|\beta,T)\ge h(\U|\beta,T)$ and
$h_\mu^-(\U|\beta,T)=h(\U|\beta,T)$.
\end{theo}
{\bf Proof:} First we prove in the 0-dimensional case that for every $\U=\{U_1,\dots,U_d\}\in \C_X$ and $\beta\in\P_\mu$ there exists $\mu\in\M_T(X)$ such that $h_\mu^+(\U|\beta,T)\ge h(\U|\beta,T)$.\\

We will then use the now classic technique used in \cite{BGH}, that works only for $h^-_\mu(\U|\beta,T)$ in our case. This is we prove the result first in the zero-dimensional case and then extend it to the general case.\\

The set of clopen sets
of a zero dimensional set $X$ is a countable set, thus the family of partitions in $\U^*$ consisting of clopen sets is countable and lets enumerate it as $\conjunto{\alpha_l}{l\in \NN}$.

Now fix $n\in\NN$ and use lemma (\ref{lemmappio+}) with $N=n^2$, $K=n$ and the family $\{T^{-n}\alpha_l\}_{l=1}^n$ applied to the cover $T^{-n}\U$ and the partition $T^{-n}\beta$. We obtain $C_n\in \beta_n^{n^2+n-1}$ and $P_n\subseteq C_n$ such that
$|P_n|\ge \frac{\mathcal{N}(\mathcal{U}_n^{n^2+n-1}|\beta_n^{n^2+n-1})}{n}$ and every element of $(\alpha_l)_n^{n^2+n-1}$
contains at most one point of $P_n$ for every $1\le l\le n$. To simplify the size of the formulas denote $N=n^2+n$. \\

Denote as $\nu_n$ the counting measure over $P_n$ and choose $1\le i,l\le n$. By definition and basic properties:

\begin{equation}
\label{equ0}
\begin{split}
H_{T^i\nu_n}((\alpha_l)_0^{N-1}|\beta_0^{N-1})&\ge H_{T^i\nu_n}((\alpha_l)_0^{N-1}|\beta_n^{N-1})-H_{T^i\nu_n}(\beta_0^{n-1}),\\
  &= H_{T^i\nu_n}((\alpha_l)_0^{N-1})-H_{T^i\nu_n}(\beta_0^{n-1}),\\
&\ge  H_{T^i\nu_n}((\alpha_l)_0^{N-1})-\log|\beta_0^{n-1}|.
\end{split}
\end{equation}

From proposition \ref{propNcond} parts  (\ref{propNcond4}), (\ref{propNcond3}) and
 (\ref{propNcond6}):

\begin{equation}
\label{equ1}
\begin{split}
\N(\U_0^{N-1}|\beta_0^{N-1})&\le \N(\U_0^{n-1}|\beta_0^{N-1})\cdot \N(\U_n^{N-1}|\beta_0^{N-1}),\\
&\le \N(\U_0^{n-1}|\beta_0^{n-1})\cdot \N(\U_n^{N-1}|\beta_n^{N-1}),\\
&\le  d^n \cdot \N(\U_n^{N-1}|\beta_n^{N-1}).
\end{split}
\end{equation}

 Since $T^{-i}(\alpha_l)_0^{N-1}\succeq (\alpha_l)_n^{N-1}$ then every element contains at most one atom of the discrete measure $\nu_n$ and so using equations (\ref{equ0}) and (\ref{equ1}):

\begin{equation}
\label{equ2}
\begin{split}
H_{T^i\nu_n}((\alpha_l)_0^{N-1}|\beta_0^{N-1})&\ge H_{\nu_n}(T^{-i}(\alpha_l)_0^{N-1})-\log|\beta_0^{n-1}|,\\
&\ge\log \left[\N(\U_n^{N-1}|\beta_n^{N-1})\over n\right]-\log|\beta_0^{n-1}|,\\
&\ge \log \left[\N(\U_0^{N-1}|\beta_0^{N-1})\over nd^n\right]-\log|\beta_0^{n-1}|.
\end{split}
\end{equation}

Fix $m\in \NN$ with $m\le n$ and decompose $N=km+b$ with $0\le b\le m-1$. Then:

\begin{equation*}
\begin{split}
H_{T^i\nu_n}((\alpha_l)_0^{N-1}|\beta_0^{N-1})&=H_{T^i\nu_n}\left(\bigvee\limits_{j=0}^{k-1}T^{-mj}(\alpha_l)_0^{m-1}\vee (\alpha_l)_{km}^{N-1}|\beta_0^{N-1}\right),\\
&\le\sum\limits_{j=0}^{k-1} H_{T^i\nu_n}(T^{-mj}(\alpha_l)_0^{m-1}|\beta_0^{N-1})+H_{T^i\nu_n}((\alpha_l)_{km}^{N-1}),\\
&\le\sum\limits_{j=0}^{k-1} H_{T^i\nu_n}(T^{-mj}(\alpha_l)_0^{m-1}|T^{-mj}\beta_0^{m-1})+\log|(\alpha_l)_{km}^{N-1}|,\\
&\le\sum\limits_{j=0}^{k-1} H_{T^{i+mj}\nu_n}((\alpha_l)_0^{m-1}|\beta_0^{m-1})+m\log d.
\end{split}
\end{equation*}

Adding up for every $0\le i\le m-1$ from equation (\ref{equ2}) we obtain:

\begin{equation}
\label{equ4}
\begin{split}
m\log \left[\N(\U_0^{N-1}|\beta_0^{N-1})\over nd^n\right]-m\log |\beta_0^{n-1}|&\le \sum\limits_{i=0}^{m-1}H_{T^i\nu_n}((\alpha_l)_0^{N-1}|\beta_0^{N-1}),\\
&\le \sum\limits_{i=0}^{m-1}\sum\limits_{j=0}^{k-1} H_{T^{i+mj}\nu_n}((\alpha_l)_0^{m-1}|\beta_0^{m-1})+m^2\log d,\\
&\le\sum\limits_{i=0}^{N-1}H_{T^i\nu_n}((\alpha_l)_0^{m-1}|\beta_0^{m-1})+m^2\log d.
\end{split}
\end{equation}

Denote as $\mu_n$the Cesaro Mean Measure of $\{T^i\nu_n\}_{i=0}^{N-1}$. By concavity and using equation (\ref{equ4}) we get for any $0\le l\le n$:

\begin{align}
\label{equ5}
\begin{split}
H_{\mu_n}((\alpha_l)_0^{m-1}|\beta_0^{m-1})&\ge {1\over N}\sum\limits_{i=0}^{N-1}H_{T^i\nu_n}((\alpha_l)_0^{m-1}|\beta_0^{m-1}),\\
&\ge{m\over N}\left(\log \left[\N(\U_0^{N-1}|\beta_0^{N-1})\over nd^n\right]-m\log d-\log|\beta_0^{n-1}|\right).
\end{split}
\end{align}

We can assume by taking a subsequence that there exists $\mu\in\M(X)$ such that in the weak-* topology $\lim\limits_{n\to\infty}\mu_n=\mu$.
Also by construction $\mu\in\M_T(X)$. Clearly:

\begin{equation}
\label{equ6}
h(\U|\beta,T)=\lim\limits_{n\to\infty}{m\over N}\left(\log \left[\N(\U_0^{N-1}|\beta_0^{N-1})\over nd^n\right]-m\log d-\log|\beta_0^{n-1}|\right).
\end{equation}

Since all elements of $(\alpha_l)_0^{m-1}$ are clopen by upper semicontinuity and equation (\ref{equ5}) divided by $m$ we obtain for any $l,m\in\NN$:

\begin{equation}
\label{equ7}
{1\over m}H_{\mu}((\alpha_l)_0^{m-1}|\beta_0^{m-1})\ge\limsup\limits_{n\to\infty}{1\over m}H_{\mu_n}((\alpha_l)_0^{m-1}|\beta_0^{m-1})= h(\U|\beta,T).
\end{equation}

Taking limit as $m$ tends to infinity proves that $h_\mu(\alpha_l|\beta,T)\ge h(\U|\beta,T)$. Since $\{\alpha_l\}_{l\in\NN}$ is dense in $L^1_\mu$ in the set of borel partitions finer than $\U$ this proves that $h^+_\mu(\U|\beta,T)\ge h(\U|\beta,T)$.  By lemma \ref{lem:pvc-df>} part (\ref{limtmhmcond}) this proves that  $h^-_\mu(\U|\beta,T)\ge h(\U|\beta,T)$ but since for any $\mu\in\M_T(X)$, $h^-_\mu(\U|\beta,T)\le h(\U|\beta,T)$ we conclude that
$h^-_\mu(\U|\beta,T)=h(\U|\beta,T)$.\\

As in \cite{BGH} we use the existence of a zero dimensional extension $(Y, S)$ of $(X,T)$, this is a continuous surjective map $\funcion{\varphi}{Y}{X}$ such that $\varphi\circ S = T\circ \varphi$ where $(Y,S)$ is a zero-dimensional TDS.\\

Now choose $\U\in \C_X$, $\beta\in \P_X$ and apply the previous case for $(Y,S)$, the cover $\varphi^{-1}\U$ and the partition $\varphi^{-1}\beta$ to obtain $\nu\in \M_S(Y)$ such that
$h^-_\nu(\varphi^{-1}\U|\varphi^{-1}\beta,S)=h(\varphi^{-1}\U|\varphi^{-1}\beta,S)$.\\

Define $\mu=\nu \circ \varphi$ from proposition \ref{hmufactors} we conclude that:

\begin{equation}\label{hmu+ineq}
h^-_\mu(\U|\beta,S)=h^-_\nu(\varphi^{-1}\U|\varphi^{-1}\beta,S)=h(\varphi^{-1}\U|\varphi^{-1}\beta,S)=h(\U|\beta,T).
\end{equation}

In order to prove the inequality for $h_\mu^+$ just notice that it is always true that $h_\mu^+(\U|\beta,T)\ge h_\mu^-(\U|\beta,T)$.

\qed

\section{Ergodic Decomposition and the equality of + and -}
There are several approaches to define the ergodic decomposition of the measure theoretical entropy and we choose the simple and clear approach taken in \cite{D}. The idea is to use the notion of disintegration of the measure $\mu$ with respect to a borel sub sigma-algebra $\D\subseteq \B_X$. There is a one to one correspondence between a sub sigma-algebra $\D$ and a measure theoretical $\varphi$-factor $(Y,\D',\nu)$ where the elements $y\in Y$ are the atoms of $\D$ and the map $\funcion{\varphi}{X}{Y}$ is defined by inclusion with $\varphi^{-1}\D'=\D$ and $\nu=\mu\circ \varphi^{-1}$.\\

For any $B\in \B$, the function $\funcion{\EE_\mu(\uno_B|\D)}{X}{\RR}$ for $\nu$-almost every $y\in Y$, is $\mu$-almost surely constant in $\varphi^{-1}y\in Y$ so we can consider that $\EE_\mu(\uno_B|\D)$ is defined over $Y$. Using this we define the disintegration of $\mu$ as the $\nu$-almost everywhere defined assignment $y\rightarrow \mu_y(\cdot):=\EE_\mu(\uno_{(\cdot)}|\D)(y)\in \M(X)$ where $\mu_y$ is a probability measure on $\B$ supported by the atom $y$. By simple properties of the conditional expectation:

\begin{equation}
\label{dismu}
\forall B\in \B, \mu(B)=\EE_\nu(\EE_\mu(\uno_B|\D))=\int\limits_{Y}\mu_y(B)\,d\nu(y)=\EE_\nu(\EE_\mu(\uno_B|\D)).
\end{equation}

In our case we consider $\D=\conjunto{B\in \B}{T^{-1}B=B}$. For the case of any $T$-invariant $\sigma$-algebra we can also define a $\varphi$-factor $(Y,\D,\nu,S)$ such that for any $y\in Y$, $\mu_{Sy}=\mu_{y}\circ T^{-1}$. In our choice of $\D$ we also have that for any $y\in Y$, $\mu_y\in \M_T^e(X)$.\\

From Theorem 2.6.4 in \cite{D} for any $\alpha\in \P_X$:

\begin{equation}
h_\mu(\alpha,T)=\int\limits_{Y}h_{\mu_y}(\alpha,T)\,d\nu(y).
\end{equation}

For any $\beta\in \P_X$ since $h_\mu(\alpha|\beta,T)=h_\mu(\alpha\vee\beta,T)-h_\mu(\beta,T)$ we also have that:

\begin{equation}
\label{ergdecalphabeta}
h_\mu(\alpha|\beta,T)=\int\limits_{Y}h_{\mu_y}(\alpha|\beta,T)\,d\nu(y).
\end{equation}



Using the same ideas given in \cite{HMRY} we prove that:

\begin{lemma}
\label{decerghmu+-}
For any $\U\in \C_X$, $\beta\in \P_X$

\begin{equation}
h^+_\mu(\U|\beta,T)=\int\limits_{Y}h^+_{\mu_y}(\U|\beta,T)\,d\nu(y),\>h^-_\mu(\U|\beta,T)=\int\limits_{Y}h^-_{\mu_y}(\U|\beta,T)\,d\nu(y).
\end{equation}
\end{lemma}

{\bf Proof:}
Let $\U=\{U_1,\dots,U_M\}\in \C_X$ and $\beta\in\P_X$. Since $X$ is a compact metric space there exists a sequence $\{\alpha_k\}_{k\in\NN}$ in $\U^*$ that is $L^1_\nu$ dense in $\U^*$ for every $\eta\in \M(X)$.

So in particular we have that for any $\eta\in\M_T(X)$:

\begin{equation}
\label{e1}
h_\eta^+(\U|\beta,T)=\inf\limits_{k\in\NN}h_\eta(\alpha_k|\beta,T).
\end{equation}

Denote for any $k\in\NN$, $\alpha_k=\{A^k_1,\dots,A^k_M\}$ where $A_m^k\subseteq U_m$ for any $1\le m\le M$. By equations (\ref{ergdecalphabeta}) and (\ref{e1}) and Fatou's Lemma we have that:

\begin{equation}
\begin{split}
h_\mu^+(\U|\beta,T)&=\inf\limits_{k\in\NN}h_\mu(\alpha_k|\beta,T)=\inf\limits_{k\in\NN}
\int\limits_{Y}h_{\mu_y}(\alpha|\beta,T)\,d\nu(y),\\
&\ge \int\limits_{Y}\inf\limits_{k\in\NN}h_{\mu_y}(\alpha|\beta,T)\,d\nu(y)=\int\limits_{Y}h_{\mu_y}^+(\U|\beta,T)\,d\nu(y).
\end{split}
\end{equation}

For every $\epsilon>0$ and $n\in\NN$ define $B_n^\epsilon=\conjunto{y\in Y}{h_{\mu_y}(\alpha|\beta,T)<h^+_{\mu_y}(\U|\beta,T)+\epsilon}$.
By equation (\ref{e1}) we know that $\mu(Y\Delta\bigcup\limits_{n\in\NN}B_n^\epsilon)=0$ and so there exists $\{Y_n\}_{n\in\NN}$ a $\nu$-partition of $X$, $\mu(Y_n)>0$ and a subsequence of partitions $\{\alpha_{k_n}\}_{k\in\NN}$ such that for any $n\in\NN$ and $\nu$-almost every $y\in Y_n$, $h_{\mu_y}(\alpha_{k_n}|\beta,T)<h^+_{\mu_y}(\U|\beta,T)+\epsilon$.\\

Now define for any $n\in\NN$ a measure $\mu_n\in \M_T(X)$ as:

\begin{equation}\label{defmun}
\forall B\in\B, \mu_n(B)={1\over \mu(Y_n)}\int\limits_{Y}\mu_y(B\cap Y_n)\,d\nu(y).
\end{equation}

By definition:

\begin{equation}\label{hkmun}
\begin{split}
h_{\mu_n}(\alpha_{k_n}|\beta,T)&={1\over \mu(Y_n)}\int\limits_{Y}h_{\mu_y}(\alpha_{k_n}|\beta,T)\,d\nu(y),\\
&\le {1\over \mu(Y_n)}\int\limits_{Y}h^+_{\mu_y}(\U|\beta,T)\,d\nu(y)+\epsilon.
\end{split}
\end{equation}

Notice that for any $n,m\in \NN$, $\mu_n(Y_m)=\uno_{\{n=m\}}$. For any $1\le m\le M$ define $A_m=\bigcup\limits_{n\in\NN} (Y_n\cap A_m^{k_n})$ and
$\alpha=\{A_1,\dots,A_M\}\in \U^*$. By construction:

\begin{equation}\label{hmual}
\begin{split}
h_{\mu}^+(\U|\beta,T)&\le h_{\mu}(\alpha|\beta,T)=\sum\limits_{n\in\NN}\mu(Y_n)h_{\mu_n}(\alpha,T),\\
&=\sum\limits_{n\in\NN}\mu(Y_n)h_{\mu_n}(\alpha_{k_n},T)\le \int\limits_{Y}h^+_{\mu_y}(\U|\beta,T)\,d\nu(y)+\epsilon.
\end{split}
\end{equation}

Since this is true for any $\epsilon>0$ this proves the reverse inequality for $h_\mu^+(\U|\beta,T)$. Finally using lemma \ref{lem:pvc-df>} part (\ref{limtmhmcond}) proves the result for $h_\mu^-(\U|\beta,T)$.\\
\qed

The final step in this procedure is to prove equality between $h_\mu^+$ and $h_\mu^-$ we state a universal version of Rohlin's Lemma following the ideas in \cite{GW}.\\

\begin{lemma}
\label{rohlin}
Let $(X,T)$ an invertible TDS. Then $\forall N\in \NN$ and
$\epsilon>0$, there exists $D\in \mathcal{B}(X)$
such that $D,TD,\ldots ,T^{N-1}D$ are pairwise disjoint and
$\mu(\bigcup_{n=0}^{N-1}T^nD)>1-\epsilon$ for every non atomic measure $\mu\in \mathcal{M}^e_T(X)$.
\end{lemma}

\begin{theo}
\label{varprinc+}
For every TDS $(X,T)$, $\U\in \C_X$, $\beta\in \P_X$ we have that $\inf\limits_{\alpha\succeq \U}\sup\limits_{\mu\in\M_T(X)}h_\mu(\alpha|\beta,T)\le h(\U|\beta,T)$.
\end{theo}
{\bf Proof:}

By Lemma (\ref{decerghmu+-}) we just need to prove it for ergodic measures moreover just for non atomic ergodic measures since atomic ergodic measures have zero measure theoretical entropy.\\

Let $\U=\{U_1,\ldots ,U_k\}$ and for any $N\in\NN$ let $\alpha'\succeq \U_0^{N-1}$ such that
$\mathcal{N}(\alpha'|\beta_0^{N-1})=\mathcal{N}(\U_0^{N-1}|\beta_0^{N-1})$. Every element of $\alpha'$ is of the form

\begin{equation*}
P=P_{i_0i_1\ldots i_{N-1}}\subseteq
\bigcap_{j=0}^{N-1}T^{-j}U_{i_j}.
\end{equation*}

Define the partition of $D$, $\alpha'_D=\{A'\cap D:A'\in \alpha'\}$. For any $\epsilon>0$ there exists
$N\in \NN$ big enough such that:

\begin{equation}\label{ek1}
\begin{split}
&\mathcal{N}(\U_0^{N-1}|\beta_0^{N-1})\le
e^{N(h(\U|\beta)+\epsilon)},\\
&-\frac{1}{N}\log\frac{1}{N}-\left(1-\frac{1}{N}\right)\log\left(1-\frac{1}{N}\right)<\epsilon.
\end{split}
\end{equation}

Let $0<\delta<1$ small enough so that
$2\sqrt{\delta}\log 2k<\epsilon.$

By Lemma \ref{rohlin}  $\exists D\subseteq X$ with
$D,TD,\ldots ,T^{N-1}D$ pairwise disjoint such that
$\mu(\bigcup_{i=0}^{N-1}T^iD)>1-\delta$ for every $\mu\in
\M_T^e(X)$ non atomic.\\

We use the partition $\alpha'_D$ to define a partition of cardinality at most $2k$ where we define $k$ atoms $\{A_1,\ldots ,A_k\}$ over the set
$\bigcup_{i=0}^{N-1}T^iD$ assigning to the set $A_i$ all sets
$$T^jP \quad{\rm con }\quad P=P_{i_0i_1\ldots i_j\ldots
i_{N-1}}\quad \text{where}\quad i_j=i\quad \text{and}\quad j=0,\ldots,N-1.
$$

On the rest of the space $\alpha$ we use any partition
that refines $\U$ that can be done with at most $k$ atoms. By construction, it is
not difficult to see that $\alpha_0^{N-1}\cap {D} \preceq
\alpha'_D$. So, for any $B\in \beta_0^{N-1}$:

\begin{equation}\label{ek2}
\begin{split}
|\alpha_0^{N-1}\cap D\cap B|&\le|\alpha'_D \cap B|\le |\alpha'\cap
B|,\\ &\le \mathcal{N}(\alpha'|\beta_0^{N-1})=
\mathcal{N}(\U_0^{N-1}|\beta_0^{N-1}).
\end{split}
\end{equation}

Let us fix  $\mu\in \M_T^e(X)$. We will show that

\begin{equation}\label{claim1}
  h_\mu(\alpha|\beta,T)\le h(\U|\beta,T)+3\epsilon.
\end{equation}

Let $E=\bigcup\limits_{i=0}^{N-1}T^i D$ such that
$\mu(E)>1-\delta.$ For $n> N$, define:

\begin{equation*}
G_n=\conjunto{x\in X}{\frac{1}{n}\sum_{i=0}^{n-1}\uno_E(T^ix)>1-\sqrt{\delta}}.
\end{equation*}

By definition $\int\limits_X \frac{1}{n}\sum_{i=0}^{n-1}\uno_E(T^ix)\,d\mu(x)>1-\delta$
 and since $\frac{1}{n}\sum_{i=0}^{n-1}\uno_E(T^ix)\le 1$ then:

\begin{equation*}
\mu(G_n)+(1-\sqrt{\delta})(1-\mu(G_n))\ge \int_X
\frac{1}{n}\sum_{i=0}^{n-1}\uno_E(T^ix)d\mu(x)>1-\delta\Rightarrow \mu(G_n)>1-\sqrt{\delta}.
\end{equation*}

For $x\in G_n$, define $S_n(x)=\conjunto{i\in \{0,1,\ldots
,n-1\}}{T^ix\in D}$. By definition for any $i=0,1,\ldots
,n-1$, $T^i x\notin E$ if and only if $\forall j\in \{0,1,\ldots
,N-1\}$, $i\notin S_n(x)+j$ thus:

\begin{equation}\label{ek4} \left|\{0,1,\ldots
,n-1\}\setminus \bigcup_{j=0}^{N-1}(S_n(x)+j)\right|\le n\sqrt{\delta}.
\end{equation}

Let $\F_n=\conjunto{S_n(x)}{x\in G_n}$. By disjointness of the $T$-iterations of $D$ for any $F=\{s_1,\ldots ,s_l\}\in \F_n, F\cap
(F+i)=\emptyset, i=1,\ldots ,N-1,$ so $|F|\le
\frac{n}{N}+1$. Define
$a_n=\left[\frac{n}{N}\right]+1$ then:
$$|\F_n|\le \sum_{j=1}^{a_n} {n\choose j}\le a_n {n\choose a_n}\le
n {n\choose a_n}.$$

By Stirling's formula and equation(\ref{ek1}):

\begin{equation*}
  \lim_{n\rightarrow\infty}\frac{1}{n}\log\left(n {n\choose a_n}\right)
=-\frac{1}{N}\log\frac{1}{N}-\left(1-\frac{1}{N}\right)\log\left(1-\frac{1}{N}\right)<\epsilon.
\end{equation*}

This implies that:

\begin{equation}\label{ek5}
\limsup_{n\rightarrow\infty}\frac{1}{n}\log|\F_n|\le\lim_{n\rightarrow\infty}\frac{1}{n}\log\left(n
{n\choose a_n}\right)\le \epsilon.
\end{equation}

For $F\in \F_n$ define $B_F=\conjunto{x\in G_n}{S_n(x)=F}$ and
$\gamma=\conjunto{B_F}{F\in \F_n}$ a measurable partition of $G_n$.\\

Let $H_F=\{0,1,\ldots
,n-1\}\setminus \bigcup_{j=0}^{N-1} (F+j)$ by (\ref{ek4}), $|H_F|\le
n\sqrt{\delta}$.\\

For any $B=\bigcap\limits_{i=0}^{n-1}T^{-i}B_{i}\in \beta_0^{n-1}$ define $\forall j=1,\ldots
,l-1$:

\begin{equation*}
  B^{s_j}=\bigcap_{t=0}^{N-1}T^{-t}B_{s_j+t}\quad \text{ and } B^{s_l}=\bigcap_{t=0}^{n-1-s_l}T^{-t}B_{s_l+t}
\end{equation*}

Then $B=\bigcap_{j=1}^{l-1}T^{-s_j}B(s_j)\cap T^{-s_l}B(s_l)\cap \bigcap_{r\in H_F}T^{-r}B_{r}$.\\

Since for all $j=1,\ldots ,l-1,$ $B(s_j)\in \beta_j\subseteq \beta_0^{n-1}$. Since $B_F\subseteq \bigcap_{j=1}^l T^{-s_j}D$ and equation (\ref{ek2}):

\begin{equation*}
\begin{split}
\left|\alpha_0^{n-1}\cap B_F \cap B\right| &\le|\bigvee_{j=1}^{l}
T^{-s_j}(\alpha_0^{N-1}\cap D \cap B^{s_j})\vee \bigvee_{r\in
H_F}T^{-r}(\alpha \cap B_{r})|,\\
& \le \prod_{j=1}^{l} |\alpha_0^{N-1}\cap D \cap
B^{s_j}|\cdot \prod_{r\in
H_F}|\alpha \cap B_{r}|,\\
& \le \mathcal{N}(\U_0^{N-1}|\beta_0^{N-1})^{l-1} |\alpha|^{|H_F|+N}\le \mathcal{N}(\U_0^{N-1}|\beta_0^{N-1})^{\frac{n}{N}}(2k)^{n\sqrt{\delta}+N}.
\end{split}
\end{equation*}

Then, adding over $F\in \F_n$, we have for every $B\in \beta_0^{n-1}$

\begin{equation}
\label{ek6}
\sum_{F\in \F_n}|\alpha_0^{n-1}\cap B_F \cap B|\le
|\F_n|\mathcal{N}(\U_0^{N-1}|\beta_0^{N-1})^{\frac{n}{N}}(2k)^{n\sqrt{\delta}+N}.
\end{equation}
Then by equation (\ref{ek6}):
\begin{equation}\label{ek7}
\begin{split}
H_{\mu_B}(\alpha_0^{n-1}\vee \gamma)&\le \sum_{F\in
\F_n}\sum_{C\in \alpha_0^{n-1}\cap
B_F}\phi(\mu_B(C))+\phi(\mu_B(X\setminus
G_n)),\\
&\le \log(\sum_{F\in \F_n}|\alpha_0^{n-1}\cap B_F \cap
B|+1),\\
&\le \log
(|\F_n|\mathcal{N}(\U_0^{N-1}|\beta_0^{N-1})^{\frac{n}{N}}(2k)^{n\sqrt{\delta}+N}
+1).
\end{split}
\end{equation}

A simple calculation shows that:
\begin{equation}\label{ek8}
\begin{split}
\mu(B)H_{\mu_B}(\alpha_0^{n-1}\cap \{X\setminus G_n\})
&=\mu((X\setminus G_n)\cap B)H_{\mu_{((X\setminus G_n)\cap B)}}(\alpha_0^{n-1})+\mu(B)H_{\mu_B}(\{X\setminus G_n\}),\\
&\le \mu((X\setminus G_n)\cap B)\log (2k)^n + \mu(B)\log 2.
\end{split}
\end{equation}

By equations (\ref{ek7}) and (\ref{ek8}):

\begin{equation}
\label{ek9}
\begin{split}
&H_\mu(\alpha_0^{n-1}|\beta_0^{n-1})\le H_\mu(\alpha_0^{n-1}\vee
(\gamma \cup \{X\setminus G_n\})|\beta_0^{n-1}),\\
&= \sum_{B\in\beta_0^{n-1}}\mu(B)\left[H_{\mu_B}(\alpha_0^{n-1}\vee
\gamma)+H_{\mu_B}(\alpha_0^{n-1}\cap \{{X\setminus G_n}\})\right],\\
&\le \log(|\F_n|\mathcal{N}(\U_0^{N-1}|\beta_0^{N-1})^{\frac{n}{N}}(2k)^{n\sqrt{\delta}+N}
+1)+\sqrt{\delta}\log (2k)^n +\log 2.
\end{split}
\end{equation}

By the form in which we choose  $\delta$ and equations (\ref{ek5}) and
(\ref{ek9}):

\begin{equation*}
\begin{split}
&h_\mu(\alpha|\beta,T) = \lim_{n \rightarrow
\infty} \frac{1}{n}H_\mu(\alpha_0^{n-1}|\beta_0^{n-1}),\\
&\le \limsup_{n\rightarrow \infty}\frac{1}{n}\left(\log|\F_n|+
\frac{n}{N}\log \mathcal{N}(\U_0^{N-1}|\beta_0^{N-1})+
(2n\sqrt{\delta}+N)\log2k +\log 2\right),\\
&\le \frac{1}{N}\log \mathcal{N}(\U_0^{N-1}|\beta_0^{N-1})+2\epsilon.
\end{split}
\end{equation*}

Finally by equation (\ref{ek1}), $h_{\mu}(\alpha|\beta,T)\le h(\U|\beta,T)+3\epsilon$.\\

Since this is true $\forall\mu\in \M_T^e(X)$ and $\epsilon$ is arbitrary we conclude that:

\begin{equation*}
\inf\limits_{\alpha\succeq \U}\sup\limits_{\mu\in \M_T^e(X)}h_{\mu}(\alpha|\beta,T)\le h(\U|\beta,T).
\end{equation*}

\qed

Fix $\beta\in\P_X$ and $\U\in \C_X$ as in \cite{GW} define:

\begin{align}
\check{h}(\U|\beta,T)&=\sup\limits_{\mu\in\M_T(X)}\inf\limits_{\alpha\succeq \U}h_\mu(\alpha|\beta,T)=\sup\limits_{\mu\in\M_T(X)}h^+_\mu(\U|\beta,T),\\
\hat{h}(\U|\beta,T)&=\inf\limits_{\alpha\succeq \U}\sup\limits_{\mu\in\M_T(X)}h_\mu(\alpha|\beta,T).
\end{align}

Clearly $\check{h}(\U|\beta,T)\le \hat{h}(\U|\beta,T)$. By Theorem \ref{varprinc} $h(\U|\beta,T)\le \check{h}(\U|\beta,T)$ and by Theorem \ref{varprinc+} $\hat{h}(\U|\beta,T)\le h(\U|\beta,T)$. So $h(\U|\beta,T)=\check{h}(\U|\beta,T)=\hat{h}(\U|\beta,T)$.\\

We can state an analogue of lemma 9 in \cite{HMRY}:

\begin{lemma}
\label{conth+}
For any $\beta\in\P_X$, for every $M\in\NN$ and $\epsilon>0$ there exists $\delta>0$ such that for every pair of measurable covers
$\U=\{U_1,\dots,U_M\}$ and $\V=\{V_1,\dots,V_M\}$ such that $\mu(\U\Delta \V)<\delta$ one has that $|h_\mu^+(\U|\beta,T)-h_\mu^+(\V|\beta,T)|\le\epsilon$.
\end{lemma}

{\bf Proof:}
Fix $M\in\NN$ and $\epsilon>0$ and choose $\delta>0$ given by lemma \ref{conthcond}.
Now fix two measurable covers
$\U=\{U_1,\dots,U_M\}$ and $\V=\{V_1,\dots,V_M\}$ such that $\mu(\U\Delta \V)<\delta$.\\

First we prove that $\alpha=\{A_1,\dots,A_M\}\in \U^*$ there exists $\alpha'\in\V^*$ such that:

\begin{equation}
\label{claim2}
h_\mu(\alpha|\beta,T)\ge h_\mu(\alpha'|\beta,T)-\epsilon.
\end{equation}

Define $\alpha'=\{A_1',\dots,A_M'\}$ as:
\begin{equation*}
\begin{split}
   A_1'&=V_1\backslash\left(\bigcup\limits_{m>1}A_m\cap V_m\right),\\
   A_m'&=V_m\backslash\left(\bigcup\limits_{k>m}A_k\cap V_k\cup\bigcup\limits_{l<m}A_l'. \right)\hbox{ for $1<m\le M$.}
\end{split}
\end{equation*}

For $1\le m\le M$ by construction $A_m\cap V_m \subseteq A_m'$, $A_m\backslash A_m'\subseteq U_m\backslash V_m$ and  $A_m'\backslash A_m\subseteq \bigcup\limits_{k\ne m}U_k\backslash V_k$. This implies that $A_m\Delta A_m'\subseteq \bigcup\limits_{m=1}^M U_m\Delta V_m$ and so $\mu(\alpha\Delta\beta)<M\cdot \mu(\U\Delta \V)<\delta$ that proves $h_\mu(\alpha|\beta,T)\ge h_\mu(\alpha'|\beta,T)-\epsilon$.\\

Since for every $\alpha'\succeq \V$, $h_\mu(\alpha'|\beta,T)\ge h_\mu^+(\V|\beta,T)$ from equation (\ref{claim2}) we conclude that $h_\mu(\alpha|\beta,T)\ge h_\mu^+(\V|\beta,T)-\epsilon$. Taking infima over all $\alpha\in \U^*$ from we conclude that
$h^+_\mu(\U|\beta,T)\ge h_\mu^+(\V|\beta,T)-\epsilon$.\\

Exchanging the roles of $\U$ and $\V$ shows that $h^+_\mu(\V|\beta,T)\ge h_\mu^+(\U|\beta,T)-\epsilon$ and so
$|h^+_\mu(\U|\beta,T)-h_\mu^+(\V|\beta,T)|\le \epsilon$.\\
\qed



\begin{theo}
\label{general h+=h-}
For every TDS $(X,T)$, $\U\in \C_X$, $\beta\in \P_X$ and $\forall\mu\in\M_T(X)$, $h_\mu^+(\U|\beta,T)=h_\mu^-(\U|\beta,T)$.
\end{theo}

{\bf Proof:}

First consider the case when $(X,T)$ is uniquely ergodic. Fix $\V\in \C_X$ and define $R=|\V|$ by lemma \ref{conth+} there exists $\delta_1>0$ such that for every $\U_1\in \C_X$, $|\U_1|=R$ and $\mu(\U_1\Delta \V)<\delta_1$ one has:

\begin{equation}
h_\mu^+(\V|\beta,T)\le h_\mu^+(\U_1|\beta,T)+\epsilon.
\end{equation}

By lemma \ref{lem:pvc-df>} part (\ref{limtmhmcond}) we can choose $M\in\NN$ such that:

\begin{equation}
\label{ec1}
{1\over M}h_\mu^+(\V_0^{M-1}|\beta,T^M)\le h_\mu^-(\U|\beta,T)+{\epsilon\over 2}.
\end{equation}

Once again by lemma \ref{conth+} there exists $\delta_2>0$ such that for every $\U_2\in \C_X$, $|\U_2|=R^M$ and $\mu(\U_2\Delta \V_0^{M-1})<\delta_1$ one has:

\begin{equation}
\label{ec2}
{1\over M}h_\mu^+(\V_0^{M-1}|\beta,T^M)\le {1\over M}h_\mu^+(\U_2|\beta,T^M)+{\epsilon\over 2}.
\end{equation}

Now choose any $\U\in \C^o_X$ with $|\U|=R$ by unique ergodicity and Theorems \ref{varprinc} and \ref{varprinc+} we have that:

\begin{equation}\label{ec3}
h_\mu^+(\U|\beta,T)=h(\U|\beta,T)=h_\mu^-(\U|\beta,T).
\end{equation}

Choose $\delta<\min\{\delta_1,{\delta_2\over R^M}\}$. By simple calculations:

\begin{equation}
\mu(\U_0^{M-1}\Delta \V_0^{M-1})\le R^M \mu(\U\Delta \V)<\delta_2.
\end{equation}

By equations (\ref{ec1}) and (\ref{ec2}):

\begin{equation}
\label{ec4}
{1\over M}h_\mu^+(\U_0^{M-1}|\beta,T^M)\le {1\over M}h_\mu^+(\V_0^{M-1}|\beta,T^M)+{\epsilon\over 2}\le h_\mu^-(\V|\beta,T)+\epsilon.
\end{equation}


By equations (\ref{ec1}), (\ref{ec2}), (\ref{ec3}) and lemma lemma \ref{lem:pvc-df>} part (\ref{igmucondtm}):

\begin{equation}
\begin{split}
h_\mu^+(\V|\beta,T)&\le h_\mu^+(\U|\beta,T)+\epsilon=h_\mu^-(\U|\beta,T)+\epsilon={1\over M}h_\mu^-(\U_0^{M-1}|\beta,T^M)+\epsilon,\\
&\le {1\over M}h_\mu^+(\U_0^{M-1}|\beta,T^M)+\epsilon\le h_\mu^-(\V|\beta,T)+2\epsilon.
\end{split}
\end{equation}

Since this is for any $\epsilon>0$ this proves that $h_\mu^+(\V|\beta,T)=h_\mu^-(\V|\beta,T)$.\\

Now by lemma \ref{decerghmu+-} we just need to prove the equality for $\mu\in\M^e_T(X)$.\\

By the Jewett-Krieger theorem (see \cite{W}) there exists a $\varphi$-measure theoretical isomorphism with a uniquely ergodic TDS $(Y,S)$.
In general $h_\nu^+(\varphi^{-1}\U|\varphi^{-1}\beta,S)\le h_\mu^+(\U|\beta,T)$ but since $\varphi$ is an isomorphism it is also true that
$h_\nu^+(\varphi^{-1}\U|\varphi^{-1}\beta,S)\ge h_\mu^+(\U|\beta,T)$. By lemma \ref{hmufactors} $h_\mu^-(\U|\beta,T)=h_\nu^-(\varphi^{-1}\U|\varphi^{-1}\beta,S)=
h_\nu^+(\varphi^{-1}\U|\varphi^{-1}\beta,S)=h_\mu^+(\U|\beta,T)$.\\

\qed

\section{Relation with other Variational Principles}

This obviously extendes the work done in \cite{R} and \cite{GW} just by considering the case $\beta=\{X\}$. Another direct application is the work done in \cite{HYZ} where they consider a fixed $\varphi$-factor $(Y,\D,S)$ and denote $\B_\varphi=\varphi^{-1}\B_X$ to define for any $\U\in \C_X$:

\begin{equation*}
\begin{split}
H_\mu(\U|Y)&:=H_\mu(\U|\B_\varphi),\>
\N(\U|Y):=\sup\limits_{y\in Y}\N(\U\cap \varphi^{-1}\{y\}),\\
h_\mu^-(\U|Y,T)&:=\lim\limits_{N\to\infty}{1\over N}H_\mu(\U_0^{N-1}|\B_\varphi), h_\mu^+(\U|Y,T):=\inf\limits_{\alpha\succeq \U}h_\mu(\alpha|\B_\varphi,T),\\
h(\U|Y,T)&:=\lim\limits_{N\to\infty}{1\over N}\log \N(\U_0^{N-1}|\B_\varphi).
\end{split}
\end{equation*}

By some classic and simple calculations (see REMARK 4 in \cite{DS} for instance) $H_\mu(\U|\B_\varphi)=\inf\limits_{\beta\in \P_X}H_\mu(\U|\varphi^{-1}\beta)$ and
$\N(\U|Y)=\inf\limits_{\beta\in \P_X}\N(\U|\varphi^{-1}\beta)$. \\

So Theorems \ref{varprinc}, \ref{varprinc+} and \ref{general h+=h-} extend to this notion as proven in \cite{HYZ}.\\

\section{Conclusion and Final Remarks}

The results of this paper show that the variational principle once again is more local than was previously stated, now with respect to the conditioning variable. This takes us one step closer to the ultimate local variational principle, this is, for an open cover conditioned with respect to a fixed open cover. Once again, the technique developed in \cite{R} to define a ``+" and ``-" definitions of measure theoretical entropy and proving their equality remains crucial and unavoidable and works for all measurable covers.\\

The requirement of the cover to be open is once more proven to be necessary for the existence of a local variational principle. This fails even for closed covers as the pioneers of this ideas knew from the beginning for the more global principles. This shows that this is not a merely combinatorial result and the topological structure is crucial to link measure and topology. In the local case, an equality is attained for a specific invariant measure and not just an equality for the suprema over all invariant measures as for the global principles but little is known about that measure. However in many cases for the global variational principles much is known about the measure that attains the equality when it exists.\\

The extension of this result for more general actions requires extra work but at least for the case of countable discrete amenable groups by using the tool of F{\o}lner sequences it seems clear how to do it.\\

\subsection{Acknowledgment}
 The author would like to thank Daniel Coronel for a lot of work and some of the main ideas in this work and Fran\c{c}ois Blanchard, Tomasz Downarovicz, Alejandro Maass and Karl Petersen for valuable remarks and discussions. The author acknowledges the support of
Programa Basal PFB 03, CMM, Universidad de Chile.


\begin{thebibliography}{MMM}

\bibitem[AKM]{AKM} R. L. Adler, A. G. Konheim and M. H. McAndrew. Topological entropy. Trans. Amer. Math. Soc.
114 (1965), 309-319.

\bibitem[B]{B} F. Blanchard . Fully positive topological entropy and topological mixing. Symbolic Dynamics and its Applications, 135 (Contemporary Mathematics). American Mathematical Society, Providence, RI, 1992, pp. 95-105.
    
\bibitem[BL]{BL} F. Blanchard and Y. Lacroix . Zero-entropy factors of topological flows. Proc. Amer. Math. Soc. 119 (1993), 985-992.

\bibitem[BGH]{BGH} F. Blanchard, E. Glasner, B. Host, A variation on the
variational principle and applications to entropy pairs,
Ergodic Th. Dynam. Sys. 17 (1997) 29-43.

\bibitem[D]{D}T. Downarovicz, Entropy in Dynamical Systems, New Mathematical Monographs, Cambridge University, Vol 18.

\bibitem[DS]{DS} T. Downarowicz, J. Serafin, Fiber entropy and conditional variational principles in compact non-metrizable spaces, Fund. Math., 172 (3) (2002), pp. 217-247

\bibitem[DZ]{DZ}Dooley, A., Zhang, G. (2015). Local entropy theory of a random dynamical system (Vol. 233, No. 1099). American Mathematical Society.

\bibitem[G]{G} T. N. T. Goodman. Relating topological entropy and measure entropy. Bull. London Math. Soc.
3 (1971), 176-180.

\bibitem[Gw]{Gw}L. W. Goodwyn. Topological entropy bounds measure-theoretic entropy. Proc. Amer. Math. Soc.
23 (1969), 679-688.

\bibitem[Gw2]{Gw2}L. W. Goodwyn, Some counter-examples in topological entropy. Topology 11 (1972), 377-385.

\bibitem[GW]{GW}E. Glasner and B. Weiss, On the interplay between measurable and topological dynamics, in
Handbook of Dynamical Systems, Vol. 1 B(eds. B. Hasselblatt and A. Katok), Elsevier B.
V., Amsterdam, 2006, 597-648.

\bibitem [HMRY]{HMRY} W. Huang, A. Maass, P. P. Romagnoli, and X. Ye, Entropy pairs and a local
Abramov formula for a measure theoretical entropy of open covers.  Ergodic Th. Dynam. Sys. 24 (2004) 1127-1153.

\bibitem [HY]{HY} W. Huang and Y. F. Yi, A local variational principle for pressure and its applications to
equilibrium states, Israel Journal of Mathematics, 161 (2007), 29-74.

\bibitem [HYZ]{HYZ} W. Huang, X. Ye, and G. Zhang,  A local variational principle for
conditional entropy.Ergodic Th. Dynam. Sys. Volume 26(2006), 219-245.

\bibitem [HYZ2]{HYZ2} W. Huang, X. Ye, G.H. Zhang, Local entropy theory for a countable discrete amenable group action, J. Funct Anal. 261 (4) (2011) 1028-1082.

\bibitem [K]{K} A. N. Kolmogorov. A new metric invariant of transient dynamical systems and automorphisms
of Lebesgue spaces. Dokl. Akad. Sci. SSSR 119 (1958), 861-864 (Russian).

\bibitem[L]{L} F. Ledrappier, A variational principle for the topological conditional entropy, Ergodic theory
(Proc. Conf., Math. Forschungsinst., Oberwolfach, 1978), Lecture Notes in Math., vol. 729,
Springer, Berlin, 1979, pp. 78-88. MR550412 (80j:54011).


\bibitem[M]{M}{M. Misiurewicz}, Topological conditional entropy,
 Studia Math. 55 (1976) 176-200.

\bibitem[R]{R}{P.P Romagnoli}, A local variational principle for the topological entropy, Ergodic Th. and Dynam.
Sys. 23 (2003) 1601-1610.

\bibitem[S]{S} {U. Shapira}, Isr. J. Math. (2007) 158: 225. doi:10.1007/s11856-007-0012-z

\bibitem[W]{W} {P. Walters}, 1982. An Introduction to Ergodic Theory. Graduate Texts in Mathematics, vol. 79. New York: Springer-Verlag.

\bibitem[Z]{Z} Zhang, Guohua., Local variational principle concerning entropy of a sofic group action, Journal of Functional Analysis 262.4 (2012): 1954-1985.
\end{thebibliography}
\end{document}